\documentclass[11pt,twoside]{article}
\usepackage{amssymb}

\newcommand{\Z}{\mathbb{Z}}

\newcommand{\C}{\mathbb{C}}
\newcommand{\F}{\mathbb{F}}

\newtheorem{Thm}{Theorem}[section]
\newtheorem{Cor}[Thm]{Corollary}

\newtheorem{Prop}[Thm]{Proposition}

\title{The character values of Iwahori--Hecke algebras on Coxeter basis 
elements}

\author{{\sc Meinolf Geck } \\[0.5cm]
(unpublished notes from 1992-1993)
}

\date{}

\begin{document}
\maketitle

Let $W$ be a finite  Weyl group, $S \subset W$ a corresponding set  of
simple reflections, and $H$ the associated generic Iwahori--Hecke algebra over
the ring $A={\Z}[u^{1/2},u^{-1/2}]$ of Laurent polynomials in an indeterminate
$u^{1/2}$ over $\Z$. We assume that all parameters of $H$ are equal to $u$, 
i.e., the basis elements $T_w$, $w \in W$, satisfy the following relations. 
\[\begin{array}{ccll} 
T_{w} T_{w'} & = & T_{ww'}               & \mbox{\ if } l(ww')=l(w)+l(w'), \\ 
T_s^2        & = & u T_1 + (u-1) T_s & \mbox{\ for } s \in S
\end{array}\] 
where $l(w)$ is the usual length function on  $W$. Let $K$ be the field of 
fractions of $A$. It is well--known that the $K$--algebra 
$H_K=K \otimes_A H$ is semisimple and split over $K$. In \cite{GP1},
G. Pfeiffer and the author have defined {\em the character table}
of $H_K$. It now seems desirable to determine these tables in all cases.
Complete results are known, up to now, for $W$ of type $G_2$, $F_4$, $E_6$, 
$E_7$ (see \cite{habil}) and $A_n$, $n \geq 1$ (see \cite{Ram}). 
For the exceptional groups, these results were achieved by ad--hoc methods
and lengthy computer calculations.

The main problem is the computation of the character
values on elements $T_w$ where $w$ is not contained in any proper parabolic
subgroup (one can use inductive arguments for the remaining elements).
There is one distinguished class of elements of this kind, namely the
class of {\em Coxeter elements}. It is the purpose of this paper to show 
that the value of an irreducible character on $T_w$, where $w \in W$ is a 
Coxeter element, is zero or, up to sign, a power of $u^{1/2}$. We also 
determine the precise values in all cases. This is achieved by reformulating 
F. Digne's and J. Michel's work \cite{DM1} on ${\cal L}$--functions of 
Deligne--Lusztig varieties  and Shintani descent (see Section~1), and 
combining it with G. Lusztig's results \cite{LFrob} on Coxeter orbits and 
eigenspaces of Frobenius (see Section~2). 

\section{Character values and $\ell$--adic cohomology}
Let $\phi$ be an irreducible character of $H_K$ and $w \in W$. Then
$\phi(T_w)$ is a polynomial in $u^{1/2}$ with integral coefficients (see
\cite{DM1}, Section~II.2). If we specialize $u$ to a prime power $q$ it is 
possible to relate the specialized characters with ${\cal L}$--functions of 
Deligne--Lusztig varieties and Shintani descent for an associated reductive 
group defined over a finite field. In this section, we make this relation 
explicit by reformulating the main results of \cite{DM1} and \cite{LBuch}.
We assume from now on that the Weyl group $W$ is irreducible.

\noindent {\bf (1.1)} Let $G$ be a simple algebraic group, 
defined over the finite field with $q$
elements, and $F:G \rightarrow G$ be the Frobenius map corresponding to
a split ${\F}_q$--rational structure on $G$. Let $B \subset G$ be an
$F$--stable Borel subgroup of $G$ and $T \subset B$ be an $F$--stable 
maximal torus of $G$ contained in $B$. Then we may identify $W$ with 
$N_G(T)/T$. The fact that $G$ is split implies that $F$ acts trivially 
on $W$. 

Let ${\bf 1}_{B^F}^{G^F}$ be the permutation representation of $G^F$ on
the cosets of $B^F$. Fix a square root $q^{1/2}$ and let 
$f:A \rightarrow \C$, $u^{1/2} \mapsto q^{1/2}$, be the corresponding
specialization. If we consider $\C$ as an $A$--module via $f$ we obtain
an isomorphism of $\C$--algebras
\[ {\rm End}_{G^F}({\bf 1}_{B^F}^{G^F}) \cong \C \otimes_A H.\]
This in turn induces a bijection, $\phi \mapsto \chi_{\phi}$, between the 
irreducible characters of $H_K$ and the irreducible characters of $G^F$ 
appearing with non--zero multiplicity in ${\bf 1}_{B^F}^{G^F}$.

We may always assume that $q$ is large so that the degrees of unipotent
characters and the orders of $G^F$ and $T^F$, for $F$--stable maximal tori
$T \subset G$, behave like polynomials in $q$.

\noindent {\bf (1.2)} For any non--negative integer
$m \geq 0$, we have a specialization $f_m:A \rightarrow \C$ which
maps $u^{1/2}$ to $(q^{1/2})^m$. To abbreviate notation, we shall write
$a_m$ instead of $a((q^{1/2})^m)$, for any element $a \in A$. 

Let $\phi$ be an irreducible character of $H_K$. Then $\phi(T_w) \in A$ for 
all $w \in W$, and we will write $\phi_m:=f_m \circ \phi$ for the 
specialized character of the 
corresponding specialized algebra $\C \otimes_A H$ (where $\C$ is regarded
as $A$--module via $f_m$).
In particular, if $m=0$ then the specialized algebra is isomorphic to the
group algebra of $W$ over $\C$, and the map $\phi \mapsto \phi_0$ defines
a bijection between the irreducible characters of $H_K$ and $W$.

\noindent {\bf (1.3)} Given $w \in W$, let
\[ X_w :=\{ x \in G \mid x^{-1}F(x) \in wB \}/(B \cap wBw^{-1})\]
be the corresponding {\em Deligne--Lusztig variety}.
Both the finite group $G^F$ and $F$ act on $X_w$, and these two actions
are compatible with each other. For any $g \in G^F$ we define
\[ R_w^1(g):={\rm Trace}(g \mid H_c^*(X_w))\]
where $H^*_c(X_w)$ is the $\ell$--adic cohomology with compact support of 
$X_w$ with coefficients in an algebraic closure of the field of $\ell$--adic
numbers. 

Given an irreducible character $\phi$ of $H_K$, we define
\[ R_{\phi}:=\frac{1}{|W|} \sum_{w \in W} \phi_0(T_w)R_w^1\]
to be the {\em almost character} of $G^F$ associated with $\phi$.
Two unipotent characters $\chi, \chi'$ of $G^F$ are
said to lie in the same {\em family}, if there exist unipotent characters
$\chi=\chi^0,\chi^1,\ldots,\chi^n=\chi'$ and characters $\phi^1,\ldots,\phi^n$ 
of $H_K$ such that $(\chi^{i-1},R_{\phi^i}) \neq 0$ and $(R_{\phi^i},\chi^i) 
\neq 0$, for all $i \geq 1$. 

\noindent {\bf (1.4)}  Let $\Gamma$ be a finite group. Then there is a
natural action of $\Gamma$ on the set of pairs $(x,\sigma)$ where $x \in 
\Gamma$ and $\sigma$ is an irreducible character of $C_{\Gamma}(x)$.
Denote by $[x,\sigma]$ the equivalence class corresponding to such a pair,
and by $M(\Gamma)$ the set of all classes. Given $c=[x,\sigma]$,
$d=[y,\tau] \in M(\Gamma)$, we define
\[ \{c,d\}:=\frac{1}{|C_{\Gamma}(x)|}\frac{1}{|C_{\Gamma}(y)|}
\sum_g \sigma(gyg^{-1})\bar{\tau}(g^{-1}xg)\]
where $g \in \Gamma$ runs over the elements satisfying $xgyg^{-1}=gyg^{-1}x$.
The numbers $\{c,d\}$ are called {\em non--abelian Fourier coefficients}.

Let ${\cal F}$ be a family of unipotent characters.
Then there exists a finite group $\Gamma$ and a bijection $M(\Gamma)
\rightarrow {\cal F}$, $c \mapsto \chi_c^{\cal F}$, such
that \[ (\chi_c^{\cal F},R_{\phi})=\left\{ \begin{array}{cl}
\{c,d\} & \mbox{ if } \chi_{\phi}=\chi_d^{\cal F} \\
0 & \mbox{ otherwise } \end{array} \right. \]
where $\chi_{\phi}$ is the component of the induced representation
${\bf 1}_{B^F}^{G^F}$ corresponding to the character $\phi$ (see (1.1)).

\noindent {\bf (1.5)} Let ${\cal F}$ be a family of unipotent characters
of $G^F$ and $\Gamma$ the associated finite group. Let $\chi \in {\cal F}$
correspond to $c=[x,\sigma] \in M(\Gamma)$. The family ${\cal F}$ will 
be called {\em exceptional} if it contains a character
$\chi$ appearing with non--zero multiplicity in ${\bf 1}_{B^F}^{G^F}$
which corresponds either to one of the characters of degree 512 of the
Weyl group of type $E_7$ or to one of the four characters of degree 4096 of
the Weyl group of type $E_8$.

Denote by $H_c^{*}(X_w)_{\chi}$ the $\chi$--isotypic component. By 
\cite{DM1}, Th\'{e}or\`{e}me~III, (2.3), there exists a complex root of
unity $\omega_{\chi}$ such that the eigenvalue of $F$ on $H_c^{*}(X_w)_{\chi}$
equals $\omega_{\chi}$ times a non--negative power of $q^{1/2}$.
By \cite{DM1}, Proposition~VII, (1.5), this root of unity can also be
described as follows.
\[ \omega_{\chi}=\left\{ \begin{array}{cl} \sigma(x)/\sigma(1) & \mbox{
if $x=1$ or ${\cal F}$ is not exceptional} \\ 
i \sigma(x)/\sigma(1) & \mbox{ otherwise.} \end{array} \right. \]

\addtocounter{Thm}{5}
\begin{Thm}
Let $\phi$ be an irreducible character of $H_K$ and $m \geq 1$. Then
\[\phi_m(T_w)=\sum_{\chi} {\rm Trace}(F^m \mid
H_c^{*}(X_w)_{\chi})\bar{\omega}_{\chi}^m(R_{\phi},\chi)\]
where $\chi$ runs over all unipotent characters of $G^F$ which lie in the
family determined by $R_{\phi}$.
\end{Thm}
{\em Proof.} Following \cite{DM1} we define, for any $g \in G^F$ and any 
positive integer $m \geq 1$, 
\[ N_w^m(g):={\rm Trace}(gF^m \mid H^*_c(X_w)).\]
By (1.5), \cite{DM1}, Definition~III, (2.1), and Th\'{e}or\`{e}me~III, (3.5), 
we have \[ N_w^m=\sum_{\phi'} \phi_m'(T_w)\sum_{\chi'} (R_{\phi'},\chi')
\omega_{\chi'}^m\chi'\]
where $\phi'$ runs over all irreducible characters of $H_K$ and $\chi'$ runs
over all unipotent characters of $G^F$. We take the scalar product of both 
sides with a unipotent character $\chi$ of $G^F$ and obtain:
\begin{eqnarray*} (N_w^m,\chi) & = & 
\sum_{\phi'}\phi_m'(T_w)\sum_{\chi'}(R_{\phi'},\chi')
\omega_{\chi'}^m(\chi',\chi) \\ & = & 
\sum_{\phi'}\phi_m'(T_w)(R_{\phi'},\chi) \omega_{\chi}^m. 
\end{eqnarray*}
Now we multiply both sides with $\bar{\omega}_{\chi}^m(R_{\phi},\chi)$ and
sum over all $\chi$:
\begin{eqnarray*} \sum_{\chi} (N_w^m,\chi)\bar{\omega}_{\chi}^m(R_{\phi},\chi)
& = & \sum_{\chi,\phi'}\phi_m'(T_w)(R_{\phi'},\chi)(R_{\phi},\chi) \\ &= & 
\sum_{\phi'}\phi_m'(T_w)\sum_{\chi} (R_{\phi'},\chi) (R_{\phi},\chi)
\\ & = & \sum_{\phi'} \phi_m'(T_w) \delta_{\phi',\phi} \mbox{ by
\cite{DM1}, Corollaire~III, (3.2)} \\
& = & \phi_m(T_w). \end{eqnarray*}
Now we observe that $(N_w^m,\chi)$ equals the trace of $F^m$ on
the $\chi$--isotypic component of $H_c^{*}(X_w)$, see \cite{DM1}, p.48.
The proof is complete. \hfill $\Box$

\section{Character values on Coxeter basis elements}
Choose a labelling $S=\{s_1,\ldots,s_n\}$ and let $w:=s_1 \cdots s_n \in W$.
Then $w$ is called a {\em Coxeter element} of $W$. It is well--known that
all the Coxeter elements corresponding to the different labellings of the
elements in $S$ are conjugate in $W$. 
Clearly, all Coxeter elements have minimal length in the conjugacy class of
$W$ in which they lie. Using \cite{GP1} we therefore get that, 
if $w,w' \in W$ are Coxeter elements then the corresponding basis elements
$T_w$, $T_{w'}$ of $H$ are conjugate by a unit in $H$. In particular, the
values of an irreducible character of $H_K$ on all Coxeter elements are equal.
In this section we apply Theorem~1.6 to determine explicitly these values in 
all cases. The argument is heavily based on Lusztig's results on Coxeter
orbits and eigenspaces of Frobenius \cite{LFrob}.

\noindent {\bf (2.1)} Let $w \in W$ be a Coxeter element.
By \cite{LFrob}, \S~6, each eigenspace of $F$ on $\bigoplus_i H_c^i(X_w)$
is contained in $H_c^i(X_w)$, for some $i$, and gives an irreducible
representation of $G^F$; furthermore, distinct eigenspaces give 
non--isomorphic representations. The total number of these representations
is $h$, the Coxeter number of $W$. Let us fix some notation.

We denote by $\lambda_1,\ldots,\lambda_h$ the eigenvalues of $F$ and by
$\chi_1,\ldots,\chi_h$ the unipotent characters of $G^F$ arising 
from the corresponding eigenspaces. For each $j$, there is exactly one 
$i_j$ such that $\chi_j$ appears as a constituent of $H_c^{i_j}(X_w)$, and
we have $\lambda_j=\omega_j(q^{1/2})^{m_j}$ where $\omega_j$ is the complex 
root of unity already introduced in (1.5) and $m_j$ is a non--negative 
integer. Then
\[ \chi_j(1)=(-1)^{i_j}[G^F:T_w^F]\lambda_j^{-1} \prod_{j' \neq j} 
(\lambda_j-\lambda_{j'})^{-1}\]
where $T_w$ is an $F$--stable maximal torus of $G$ obtained by twisting the
maximally--split torus $T$ with $w$ (see \cite{LFrob}, Theorem~6.1(ii)). 

For each $j$, we denote by $q^{a_j}$ the maximal power of $q$ dividing
$\chi_j(1)$ (note that this is well-defined since we assume that $q$ is large
enough). Using the lists of eigenvalues of $F$ in \cite{LFrob}, (7.3), we can
determine the number $a_j$, and it turns out that it only depends on the
absolute value of $\lambda_j$. In the following tables we list these numbers
for all types of $W$ individually. 

$A_n: \begin{array}{c|cc}  |\lambda_j| & q^i & i=0,\ldots,n \\ \hline
                            a_j & \frac{1}{2}(n-i)(n-i+1) \end{array}$\\
$B_n: \begin{array}{c|cc} |\lambda_j| & q^i & i=0,\ldots,n \\ \hline
                            a_j & (n-i)^2 \end{array}$\\
$D_n: \begin{array}{c|cccc} |\lambda_j| &1&q^n& q^i & i=1,\ldots,n-1 \\ \hline
                            a_j & n(n-1) & 0 & (n-i)(n-i-1)+1 \end{array}$\\
$G_2: \begin{array}{c|ccc} |\lambda_j| & 1 & q & q^2 \\ \hline
                            a_j & 6 & 1 & 0 \end{array}$\\
$F_4: \begin{array}{c|ccccc} |\lambda_j| & 1 & q & q^2 & q^3 & q^4 \\ \hline
                            a_j & 24 & 13 & 4 & 1 & 0 \end{array}$\\
$E_6: \begin{array}{c|ccccccc} |\lambda_j| &1&q&q^2&q^3&q^4&q^5&q^6 \\
        \hline a_j & 36 & 25 & 15 & 7 & 3 & 1 & 0 \end{array}$\\
$E_7: \begin{array}{c|ccccccccc} |\lambda_j| &1&q&q^2&q^3&q^{7/2}&
         q^4&q^5&q^6&q^7 \\
        \hline a_j & 63 & 46 & 30 & 16 & 11 & 7 & 3 & 1 & 0 \end{array}$\\
$E_8: \begin{array}{c|ccccccccccc} |\lambda_j| &1&q&q^2&q^3&q^{7/2}&
         q^4&q^{9/2}&q^5&q^6&q^7&q^8 \\
        \hline a_j & 120 & 91 & 63 & 37 & 26 & 16 & 11 & 7 & 3 & 1 & 0
        \end{array}$

\vspace*{2mm}
\noindent {\bf (2.2)} Let $\phi$ be an irreducible character of $H_K$, and
${\cal F}$ the family of unipotent characters of $G^F$ determined by 
$R_{\phi}$. By \cite{LBuch}, (4.26.3), there exists a non--negative integer
$a_{\phi}$ such that $\chi(1)=cq^{a_{\phi}}N$ where $c$ is a certain
Fourier coefficient and $N$ is an integer satisfying $N \equiv 1$ modulo~$q$.
Thus $q^{a_{\phi}}$ is the exact  power of $q$ dividing $\chi(1)$, for all
$\chi \in {\cal F}$. From \cite{LBuch}, (4.1), we conclude that $a_{\phi}$ 
can also be characterized by the property that $u^{a_{\phi}}$ is the exact 
power of $u$ dividing the generic degree polynomial of $\phi$.

Now let $j \in \{1,\ldots,h\}$, $w \in W$ be a Coxeter element, and assume 
that $\chi_j$ lies in ${\cal F}$. Then $a_{\phi}=a_j$. The results in 
(2.1) show that the absolute value of $\lambda_j$ is uniquely
determined by $a_{\phi}$. We shall denote it by $(q^{1/2})^{m_{\phi}}$.
With this notation we have, for any $m \geq 1$, 
\[ {\rm Trace}(F^m \mid H_c^{*}(X_w))_{\chi_j}=(-1)^{i_j}\omega_j^m
(q^{1/2})^{mm_j}=(-)^{i_j}\omega_j^m (q^{1/2})^{mm_{\phi}}.\]
Hence, for all such $j$, these terms have the same absolute value.
Combining this with Theorem~1.6 and varying over infinitely 
many $q$ and $m$ yields the formula: 
\[\phi(T_w)=\varepsilon_{\phi}(u^{1/2})^{m_{\phi}} \quad \mbox{ where } 
\quad \varepsilon_{\phi}:= \sum_j (-1)^{i_j}(R_{\phi},\chi_j)\]
(sum over all $j$ such that $\chi_j$ lies in the family ${\cal F}$).
Since $\phi(T_w) \in A$, we also conclude that $\varepsilon_{\phi}$ is
an integer. 

\addtocounter{Thm}{2}
\begin{Prop} Suppose that $W$ is of classical type and let $w \in W$ be a 
Coxeter element.

(i) Assume that $W$ is of type $A_n$, for $n \geq 0$.
The irreducible characters of $H_K$ are parametrized
by partitions of $n+1$ (e.g., the sign character corresponds to the 
partition $(1^{n+1})$). Suppose that $\phi$ is labelled by the 
partition $\alpha$.  Then
\[ \phi(T_w)=\left\{ \begin{array}{cl} (-1)^{n+1+k}u^{k-1} & \mbox{ if }
\alpha=(1^{n+1+k},k) \\ 0 & \mbox{ otherwise.} \end{array} \right. \] 

(ii) Assume that $W$ is of type $B_n$, for $n \geq 2$.
The irreducible characters of $H_K$ are parametrized by ordered pairs
$(\alpha,\beta)$ of partitions  such that $|\alpha|+|\beta|=n$
(e.g., the sign character corresponds
the  pair $(-,1^n)$). Suppose that $\phi$ is labelled by the pair $(\alpha,
\beta)$.  Then
\[ \phi(T_w)=\left\{ \begin{array}{cl} (-1)^{n+k+1}u^{k-1} & 
\mbox{ if } \alpha=(1^{n-k},k),\beta=(-) \\ (-1)^{n+k}u^{k} & 
\mbox{ if } \alpha=(-), \beta=(1^{n-k},k)
\\ 0 & \mbox{ otherwise.}  \end{array} \right. \]

(iii) Assume that $W$ is of type $D_n$, for $n \geq 4$.
The irreducible characters of $H_K$ are parametrized by unordered pairs
$(\alpha,\beta)$ of partitions of $n$  such that $|\alpha|+|\beta|=n$, and
there are two irreducible characters if $\alpha=\beta$
(e.g., the sign character corresponds to the  pair $(-,1^n)$). 
Suppose that $\phi$ is labelled by the pair $(\alpha, \beta)$. Then 
\[ \phi(T_w)=\left\{ \begin{array}{cl} 
(-1)^{n+k}u^k & \mbox{ if } \alpha=(1),\beta=(1^{n-k-1},k) 
\mbox{ for $k \geq 1$} \\ 
(-1)^{n+k+1}u^k &\mbox{ if } \alpha=(-),\beta=(1^{n-k-2},2,k) 
\mbox{ for $k \geq 2$} \\ 
(-1)^n & \mbox{ if } \alpha=(-),\beta=(1^n) \\
u^n & \mbox{ if } \alpha=(-),\beta=(n) \\
0 & \mbox{ otherwise.}  \end{array} \right. \]
For the labelling of the irreducible characters of $W$ by partitions or pairs 
of partitions, we refer to \cite{C2}, Section~11.4.
\end{Prop}
{\em Proof.} (a) At first we note that it will be sufficient to prove
that $\phi(T_w) \neq 0$ for the specified labels of $\phi$. This can be seen
as follows. We have $\phi(T_w)=\varepsilon_{\phi}(u^{1/2})^{m_{\phi}}$ by 
(2.2), hence specialization $u^{1/2} \mapsto 1$ shows that $\varepsilon_{\phi}$
equals the value of the specialized character $\phi_0$ of $W$ on~$w$. 
By the second
orthogonality realation for the characters of $W$, we have $\sum_{\phi}
\varepsilon_{\phi}^2=|C_W(w)|$, and this equals the Coxeter number $h$ of $W$
which is given in \cite{LFrob}, p.106. Now we check that, in each case
(i)--(iii), we have specified precisely $h$ characters for which the value
on $T_w$ is non--zero, so the remaining ones must be zero.

(b) Next we need a criterion to decide whether or not a unipotent character
of $G^F$ is a constituent of $R_w^1$. This can be done as follows.
Let $h$ be the Coxeter number of $W$ and $\Phi_h$ the $h$--th cyclotomic 
polynomial. By \cite{L1}, Lemma~3.30, and \cite{LFrob}, Corollary~6.16, we see
that a unipotent character $\chi$ is a constituent of $R_w^1$ 
if and only if $\Phi_h(q)$ does not divide the degree of $\chi$ 
(regarded as a polynomial in $q$). By using the formula for the degrees of the 
unipotent characters in \cite{C2}, Section~13.8, this can be easily checked
from the symbol labelling the character.

(c) Now assume that the family determined by $R_{\phi}$ contains only one
element. Then this must be the principal series character $\chi_{\phi}$. By
(b) we can check whether or not it is a constituent of $R_w^1$. If it is not,
then $\phi(T_w)=0$, by (2.2). So let us assume that $\chi_{\phi}=\chi_j$, for
some $j \in \{1,\ldots,h\}$. By (2.2) we then have 
$\phi(T_w)=(-1)^{i_j}(u^{1/2})^{m_{\phi}}$, since $(R_{\phi},\chi_{\phi})=1$. 
The number $m_{\phi}$ is determined
by the degree formula for $\chi_{\phi}$ and the tables in (2.1); it is
always even. So we only need to compute $i_j$. By \cite{L1}, Example~3.10(a), 
the eigenvalue of $F$ corresponding to the principal series character 
$\chi_{\phi}$ must be positive. The list of eigenvalues in \cite{LFrob}, 
(7.3), then shows that $i_j=n+k$ if $m_{\phi}=2k$.

Using these arguments we can easily settle case (i), since all unipotent 
characters lie in the principal series and all families contain just one 
element. Other proofs can be found in \cite{DM1}, Lemme~V, (3.12), or 
\cite{Ram}.

(d) The case that $W$ is of type $B_n$ can now be handled in an elementary
way, as follows. 
Let $S:=\{t,s_1,\ldots,s_{n-1}\}$ where $S'=\{s_1,\ldots,s_{n-1}\}$
generates a parabolic subgroup $W'$ of type $A_{n-1}$ and $ts_1$ has order~4.
The irreducible characters of the corresponding parabolic subalgebra $H_K'$
are parametrized by partitions of $n$. Let $\rho':H_K' \rightarrow K^{r \times
r}$ be an irreducible representation with character $\phi'$ labelled by the
partition $\alpha$. Denote by $E_r$ the $r \times r$--identity and define
$\rho^1(T_t):=uE_r$, $\rho^2(T_t):=-E_r$, $\rho^1(T_{s_i})
:=\rho^2(T_{s_i})= \rho'(T_{s_i})$ for $i=1,\ldots,n-1$.
By checking the defining relations for $H_K$ we readily see that 
$\rho^1$ and $\rho^2$
are non--isomorphic extensions of  $\rho'$ to $H_K$. Denote  by $\phi^1$ and
$\phi^2$ their characters, and let $w'=s_1 \cdots s_{n-1} \in W'$ be a 
Coxeter element. Then $w=tw'$ is a Coxeter element in $W$, and we have
\[ \phi^1(T_w)=u\phi'(T_{w'}) \mbox{ and } \phi^2(T_w)=-\phi'(T_{w'}).\]
The definition of the parametrization of the characters of $W$ shows that
$\phi^1$ is labelled by the pair $(-,\alpha)$ and $\phi^2$ is labelled by
$(\alpha,-)$. If $\alpha=(1^{n-k},k)$ for some $k \in \{1,\ldots,n\}$ then,
by (i) and the above equation, we have $\phi^1(T_w)=(-1)^{n+k}u^k$ and 
$\phi^2(T_w)=(-1)^{n+k+1}u^{k-1}$. In this way, we obtain $2n$ irreducible 
characters of $H_K$ with non--zero value on $T_w$, and we are done by (a).

(e) Now let $W$ be of type $D_n$, for $n \geq 4$. The characters are labelled
by symbols of rank $n$ with some extra conditions. This is described in 
detail in \cite{C2}, p.471ff. Given a symbol $X$, there is a unique special 
symbol $X'$ in the same family. Let $Z$ be the set of numbers which appear in 
just one row of $X'$. Then the number of unipotent characters in this family
equals $2^{|Z|-2}$. Given a character $\phi$ of $H_K$ labelled 
by a pair of partitions, the associated symbol is determined by the procedure 
described in \cite{C2}, Proposition~11.4.4.

The statement in (iii) for the sign and the index character is clear. Now let 
$\alpha=(1)$, $\beta=(1^{n-1-k},k)$ for some $k \in \{1,\ldots,n-1\}$.
The associated symbol is 
{\small \[ X_1:=\left( \begin{array}{ccccccc} 0 & 1 & 2 & \ldots & m-3 & m-2 & 
m \\ 1 & 2 & 3 & \ldots & m-2 & m-1 & n-1 \end{array} \right)\]}
where $m=n-k$ is the number of parts of $\beta$. This is a special symbol.
We have $Z=\{0,n-2\}$ if $k=1$, $Z=\{0,n-1\}$ if $k=n-1$, and 
$Z=\{0,m-1,m,n-1\}$ otherwise. Hence, if $k=1$ or $k=n-1$ then $|Z|=2$, 
and we are in the situation of (c). So let us assume, from now on, that 
$2 \leq k \leq n-2$. Then $|Z|=4$.  Consider the symbols
{\small \begin{eqnarray*} 
X_2& := & \left( \begin{array}{ccccccc} 0 & 1 & 2 & \ldots & m-3 & m-2 & m-1 \\
1 & 2 & 3 & \ldots & m-2 & m & n-1 \end{array} \right) \\
X_3& :=&\left( \begin{array}{ccccccc} 0 & 1 & 2 & \ldots & m-3 & m-2 & n-1 \\
1 & 2 & 3 & \ldots & m-2 & m-1 & m \end{array} \right) \\
X_4 & := & \left( \begin{array}{ccccccccc} 0 & 1 & 2 & \ldots & m-3 & m-2 & m-1 
& m & n-1 \\ 1 & 2 & 3 & \ldots & m-2 & & & & \end{array} \right)
\end{eqnarray*}}
Then $X_1,\ldots,X_4$ all have the same entries with the same multiplicities.
Hence they are precisely the set of symbols labelling the characters in the
family determined by $R_{\phi}$.  The corresponding Fourier transform matrix 
is given by
{\small \[\begin{array}{c|rrrr} & X_1 & X_2 & X_3 & X_4 \\  \hline
X_1 & 1/2 & 1/2 & 1/2 & 1/2 \\ X_2 & 1/2 & 1/2 & -1/2 & -1/2 \\
X_3 & 1/2 & -1/2 & 1/2 & -1/2 \\ X_4 & 1/2 & -1/2 & -1/2 & 1/2 \end{array}\]}
By using the degree formula we see that $q^k$ is the exact power of $q$ 
dividing the characters in this family. Comparison with (2.1) shows that
$m_{\phi}=2k$. Using (b) we see that the characters corresponding to the 
symbols $X_3$ and $X_4$ are constituents of $R_w^1$. The one labelled by 
$X_3$ is a principal series character, the other is not. Let $j,j'$ be 
the indices of the characters labelled by $X_3$, $X_4$, in the notation of 
(2.1). From the tables in \cite{LFrob}, (7.3), we see that $i_j=n+k$ and 
$i_{j'}=n+k-2$ (note again that a principal series character corresponds 
to a positive eigenvalue). The formula in (2.2) now reads 
\[\varepsilon_{\phi}=(-1)^{n+k}((R_{\phi},\chi_j)+(R_{\phi},\chi_{j'})).\]
The symbol associated with $\phi$ is $X_1$, hence we have
$\varepsilon_{\phi}=(-1)^{n+k}(1/2+1/2)=(-1)^{n+k}$. 

Finally, let $\alpha=(-)$, $\beta=(1^{n-k-2},2,k)$ for $k \in 
\{2,\ldots,n-2\}$.
We find that the associated symbol is $X_2$. So we can use the above
Fourier transform matrix once more and obtain $\varepsilon_{\phi}=
(-1)^{n+k}(-1/2-1/2)=(-1)^{n+k+1}$. Thus, there at least $2(n-1)$ character
of $H_K$ with non--zero value on $T_w$. Since this also is the Coxeter number
of $W$, we are done. The proof is complete. \hfill $\Box$

\vspace*{2mm}
\noindent {\bf (2.4)} Now we consider the case that $W$ is of exceptional type. 
By \cite{C2}, Section~13.2, the irreducible characters have a labelling by 
certain pairs of non--negative integers. For each case, we now give those
labels for which the corresponding irreducible character of $H_K$ 
has non--zero value on $T_w$. Note that it is important now to fix square
roots of $q$ and $1$ in order to have well-defined labellings of the
characters.

%
%
%
%

\vspace*{2mm}
\noindent $G_2$: $\begin{array}{cccccc} 
(1,0)&(1,6)&(1,3)'& (1,3)''& (2,1)& (2,2) \\ \hline
u^2&1&-u&-u&u&-u \end{array}$

\vspace*{2mm}
\noindent $F_4$: $\begin{array}{cccccc} 
(1,0)& (1,12)''&(1,12)'&(1,24)& (2,4)''&(2,16)' \\ \hline
u^4&u^2&u^2&1&-u^3&-u \end{array}$ \\
$\;\;\;\;\;\begin{array}{cccccc} 
(2,4)'&(2,16)''& (4,8)&(6,6)'&(6,6)''&(12,4) \\ \hline
-u^3&-u&u^2&-u^2&-u^2&u^2 \end{array}$

 \vspace*{2mm}
\noindent $E_6$: $\begin{array}{cccccc}
(1,0)&(1,36)&(10,9)& (6,1)&(6,25)&(20,10) \\ \hline
u^6&1&-u^3&-u^5&-u&u^3 \end{array}$\\
$\;\;\;\;\;\begin{array}{cccccc}
(15,4)&(15,16)&(30,3)& (30,15)&(60,8)&(90,8) \\ \hline
-u^4&-u^2&u^4&u^2&u^3&-u^3 \end{array}$

\vspace*{2mm}
\noindent $E_7$: $\begin{array}{cccccccc} 
(1,0)&(1,63)&(7,46)&(7,1)&(35,22)&(35,13)&(35,4)&(35,31) \\ 
\hline u^7&-1&u&-u^6&-u^3&u^4&-u^5&u^2\end{array}$
$\;\;\;\;\;\begin{array}{cccccc} 
(56,30)&(56,3)&(70,18)&(70,9)&(280,18)&(280,9) \\ \hline
-u^2&u^5&u^3&-u^4&u^3&-u^4 \end{array}$\\
$\;\;\;\;\;\begin{array}{cccc} 
(280,8)&(280,17)&(512,12)&(512,11) \\ \hline
u^4&-u^3&-u^{7/2}&u^{7/2} \end{array}$

\vspace*{2mm}
\noindent $E_8$:  $\begin{array}{cccccccc} 
(1,0)&(1,120)&(70,32)&(84,4)&(84,64)&(420,20)&(1134,20) \\ \hline
u^8&1&-u^4&-u^6&-u^2&-u^4&u^4\end{array}$\\
$\;\;\;\;\;\begin{array}{cccccc}
(1680,22)&(1344,8)&(1344,38)&(4480,16)&(4536,18)&(5670,18) \\ \hline
u^4&u^5&u^3&u^4&u^4&-u^4 \end{array}$\\
$\;\;\;\;\;\begin{array}{ccccccc}
(4096,12)&(4096,26)&(8,1)&(8,91)&(56,19)&(56,49)&(112,3) \\ \hline
-u^{9/2}&-u^{7/2}&-u^7&-u&u^5&u^3&u^6 \end{array}$\\
$\;\;\;\;\;\begin{array}{cccccc}
(112,63)&(448,25)&(448,9)&(448,39)&(1008,9)&(1008,39) \\ \hline
u^2&-u^4&-u^5&-u^3&-u^5&-u^3 \end{array}$\\
$\;\;\;\;\;\begin{array}{cccc}
(2016,19)&(7168,17)&(4096,11)&(4096,27) \\ \hline
u^4&-u^4&u^{9/2}&u^{7/2} \end{array}$

\noindent {\em Proof.} This is just a matter of going through all cases. 
The exponent $m_{\phi}$ is determined by the degree of the unipotent
characters $\chi_j$ and the tables in (2.1); the number 
$\varepsilon_{\phi}$ is the value of the specialized character $\phi_0$ of 
$W$ on $w$, so this can be read off the ordinary character  table of $W$.
For $W$ of type $G_2$ or $F_4$, this table can be found in
\cite{C2}, p.412/413; for $W$ of type $E_6$, $E_7$ or $E_8$, it is 
printed in \cite{Atl}, p.26/46/86.  

We give an example: Suppose $W$ has type $E_8$ and let $\phi$ be the 
irreducible character of $H_K$ with label $(4096,12)$. By \cite{C2}, p.485, 
we see that the integer $a_{\phi}$ defined in (2.2) is 11. By (2.1), we see 
that the absolute value of the corresponding eigenvalue of $F$ must be 
$q^{9/2}$. By (2.2), we therefore have $\phi(T_w)=
\varepsilon_{\phi}u^{9/2}$. We now specialize $u^{1/2} \mapsto 1$. Then the 
specialized character of $W$ has value $\varepsilon_{\phi}$ on $w$. 
Hence $\varepsilon_{\phi}$ is determined by the ordinary character table of 
$W$ in \cite{Atl}, p.86. 

The precise value of $\varepsilon_{\phi}$ could also be determined by using
the Fourier transform matrix of the family associated with $\phi$, exactly
in the same way as in the proof of Proposition~2.3. \hfill $\Box$

Summarizing the explicit results we have computed in Proposition~2.3 and
(2.4) we finally have:
\addtocounter{Thm}{1}
\begin{Cor}
Let $\phi$ be an irreducible character of $H_K$ and $w \in W$ a Coxeter
element. Then $\phi(T_w)=\varepsilon_{\phi}(u^{1/2})^{m_{\phi}}$ where
$m_{\phi}$ is as defined in (2.2) and $\varepsilon_{\phi} \in \{0,1,-1\}$.
\end{Cor}

\medskip
{\it Note added November 2024}. All character values of generic 
Iwahori--Hecke algebras are known since the end of the 1990s; see, e.g., 
the book by G. Pfeiffer and the author on ``Characters of finite Coxeter 
groups and Iwahori--Hecke algebras'', Oxford Univ. Press, 2000.

\end{document}